\documentclass[8pt]{article}
\title{Phase transition results for three Ramsey-like theorems \footnote{This is the pre-peer reviewed version of a paper which has been accepted for publication at the Notre Dame Journal of Formal Logic (2016)}\\ DRAFT}
 
\author{Florian Pelupessy}
\usepackage[english]{babel}
\usepackage{amsmath, amssymb}
\usepackage{url}
\usepackage{xcolor}
\usepackage{hyperref}
\usepackage{tikz}
\usepackage{multicol}
\usepackage{parskip}
\usepackage[initials]{amsrefs}
\usepackage[T1]{fontenc}
\usepackage{libertine}

\newtheorem{theorem}{Theorem}
\newtheorem{lemma}{Lemma}
\newtheorem{definition}{Definition}

\newtheorem{corollary}[theorem]{Corollary}

\newtheorem{conjecture}{Conjecture}

\begin{document}
\maketitle
\noindent
\begin{abstract}
We present a unified method for showing phase transition results for three Ramsey theorem variants.  
\end{abstract}

Phase transitions in logic are a recent development in unprovability. The general programme, started by Andreas Weiermann, is to classify parameter functions $f\colon \mathbb{N} \rightarrow \mathbb{N}$ according to the provability of a parametrised theorem $\varphi_f$ in a theory $T$. We study these transitions with the goal of gaining a better understanding of unprovability. More details on this programme, with an overview of related publications, can be found at \cite{weiermannweb}.

In this paper we will study the transition results for three Ramsey theorem variants: Friedman's finite adjacent Ramsey theorem, the Paris--Harrington theorem and the Kanamori--McAloon theorem. The latter two of these have been studied previously in \cite{weiermann2004} and \cite{CLW}, but the methods used in the present paper are a natural continuation of the method for the adjacent Ramsey theorem. The emphasis of this method is on connecting the variants $\varphi_k$ of the theorem for constant functions $k$ with the resulting classification of parameter values. Furthermore this method is independent of whether the original method of showing unprovability for the non-parametrised theory involves proof/recursion theory or model theoretic constructions. 

Additionally we will provide some general tools to streamline proofs of phase transition results: the upper bounds lemmas \ref{lemma:upperbounds}, \ref{lemma:upperboundssharpening} and the lower bounds sharpening lemmas \ref{lemma:prooftheoreticsharpening}, \ref{lemma:modeltheoreticsharpening}. The manner in which these lemmas are stated  indicates the most important steps in the proofs of (sharpened) phase transitions. This removes the need to repeat some ad-hoc arguments for each transition result, thus simplifying the proofs. 

This paper is divided into four sections. Section~\ref{section:threeramsey} introduces the three Ramsey theorem variants and the transition results. Section~\ref{section:lowerbounds} is dedicated to independence and Section~\ref{section:upperbounds} is dedicated to provability. We conclude with some observations on phase transitions in Section~\ref{section:observation}

\section{Three Ramsey-like theorems}~\label{section:threeramsey}
\subsection{Adjacent Ramsey}
The finite adjacent Ramsey theorem is one of the latest independence results at the level of $\mathrm{PA}$ and was first presented in \cite{friedman2010}. Independence of the variants with fixed dimension is examined extensively in \cite{friedmanpelupessy}. This examination uses proof theoretic techniques. Showing independence using model theoretic construcions is still an open problem. As in the case of the othere Ramsey variants we will call functions $C\colon \{0, \dots , R\}^d \rightarrow \mathbb{N}^r$ colourings. Notice the distinction between parameter functions, which are provided externally and colourings, which are quantified inside the theorems.
\begin{definition}
For $r$-tuples $a,b$:
\[
a \leq b \Leftrightarrow (a)_1 \leq (b)_1 \wedge \dots \wedge (a)_r \leq (b)_r
\]
\end{definition}
\begin{definition}
A colouring $C\colon \{0, \dots , R\}^d \rightarrow \mathbb{N}^r$ is $f$-limited if 
\[
\max C(x) \leq f(\max x)+1 \mathrm{\ for \ all \ } x \in \{0, \dots , R\}^d.
\] 
\end{definition}
\begin{theorem}[$\mathrm{AR}_f$] \label{thm:AR}
For every $d, r$ there exists $R$ such that for every limited colouring $C\colon \{0, \dots , R\}^d \rightarrow \mathbb{N}^r$ there exist $x_1 < \dots < x_{d+1} \leq R$ with $C(x_1, \dots , x_d) \leq C(x_2, \dots , x_{d+1})$. 
\end{theorem}
\begin{definition}
We denote the smallest $R$ from $\mathrm{AR}_f$ with $\mathrm{AR}_f^d(r)$. 
\end{definition}
\begin{theorem}\label{thm:independence:AR}
  ~{}
  \begin{enumerate}
  \item $\mathrm{I}\Sigma_{d+1} \nvdash \mathrm{AR}^{d+1}_{\mathrm{id}}$. 
  \item $\mathrm{PA} \nvdash \mathrm{AR}_{\mathrm{id}}$.
  \end{enumerate}
\end{theorem}
\emph{Proof:} See~\cite{friedmanpelupessy}.
\begin{flushright} $\Box$ \end{flushright}
\subsection{Paris--Harrington}
The Paris--Harrington theorem is one of the earliest examples of natural theorems which are independent of $\mathrm{PA}$.  This was first shown using model theoretic methods in~\cite{parisharrington}, later this was shown using proof theoretic methods in~\cite{ketonensolovay},~\cite{loeblnesetril} and in~\cite{friedmanpelupessy}. 
\begin{definition}
$[X]^d$ is the set of $d$-element subsets of $X$, $[m,R]^d=[\{m, \dots , R\}]^d$ and $[R]^d=[0, R]^d$.
\end{definition}
\begin{definition}
Given a colouring $C\colon [m,R]^d \rightarrow r$, we call a set $H$ homogeneous for $C$ or $C$-homogeneous if $C$ is constant on $[H]^d$.
\end{definition}
\begin{theorem}[$\mathrm{PH}_f$] \index{$\mathrm{PH}_f$} 
For every $d,r,m$ there exists an $R$ such that for every colouring $C\colon [m,R]^d \rightarrow r$ there exists an $H \subseteq [m,R]$ of size $f(\min H)$ for which $C$ limited to $[H]^d$ is constant.
\end{theorem}
\begin{definition}
We denote the smallest $R$ from $\mathrm{PH}_f$ with $\mathrm{PH}^d_f(m,r)$. We call a colouring  $C\colon [m,R]^d \rightarrow r$ bad if every $C$-homogeneous set has size strictly less than $f(\min H)$. 
\end{definition} 
\begin{theorem}\label{thm:independence:PH}
  ~{}
  \begin{enumerate}
  \item $\mathrm{I}\Sigma_{d+1} \nvdash \mathrm{PH}^{d+2}_{\mathrm{id}}$. 
  \item $\mathrm{PA} \nvdash \mathrm{PH}_{\mathrm{id}}$.
  \end{enumerate}
\end{theorem}
\emph{Proof:} See~\cite{friedmanpelupessy}.
\begin{flushright} $\Box$ \end{flushright}

\subsection{Kanamori--McAloon}
We wil examine the following variant of the Kanamori--McAloon theorem:
\begin{theorem}[$\mathrm{KM}_f$] 
For every $d,m,a$ there exists $R$ such that for every $C\colon [a, R]^d \rightarrow \mathbb{N}$ with $C(x) \leq f (\min x)$ there exists $H \subseteq R$ of size $m$ for which for all $x,y \in [H]^d$ with $\min x= \min y$ we have $C(x)=C(y)$.
\end{theorem}
\begin{definition}
We denote the smallest $R$ from $\mathrm{KM}_f$ with $\mathrm{KM}^d_f(a,m)$.
\end{definition}
\begin{theorem}\label{thm:independence:KM}
  ~{}
  \begin{enumerate}
  \item $\mathrm{I}\Sigma_{d+1} \nvdash \mathrm{KM}^{d+2}_{\mathrm{id}}$. 
  \item $\mathrm{PA} \nvdash \mathrm{KM}_{\mathrm{id}}$.
  \end{enumerate}
\end{theorem}
\emph{Proof:} See~\cite{kanamorimcaloon}.
\begin{flushright} $\Box$ \end{flushright}

\subsection{Phase transition results}
All parameter functions are assumed to be nondecreasing. For every $f\colon \mathbb{N} \rightarrow \mathbb{N}$ the inverse is:
\[
f^{-1} (i) = \min \{ j: f(j) \geq f(i)\},
\]
$\mathrm{2}_0(i)=i$, $\mathrm{2}_{n+1}(i)=2^{\mathrm{2}_n(i)}$, $\log$ is the inverse of $i \mapsto 2^i$, $\log^n$ is the inverse of $i \mapsto \mathrm{2}_n(i)$, $\sqrt[c]{\log^n}$ is the inverse of $i \mapsto \mathrm{2}_n(i^c)$ and $i \mapsto \frac{i}{c}$ is the inverse of $i \mapsto i \cdot c$, where $\frac{x}{0}=1$. We use $\varphi_f^d$ to denote the theorem $\varphi_f$ with fixed dimension $d$. 
\begin{theorem}\label{thm:transition:AR}
  ~{}
  \begin{enumerate}
  \item $\mathrm{I}\Sigma_{d+1} \nvdash \mathrm{AR}^{d+1}_{\sqrt[c]{\log^d}}$ for every $c$. 
  \item $\mathrm{PA} \nvdash \mathrm{AR}_{\log^n}$ for every $n$.
  \item  $\mathrm{I}\Sigma_{1} \vdash \mathrm{AR}^{d+1}_{\log^{d+1}}$. 
  \item $\mathrm{I}\Sigma_1 \vdash \mathrm{AR}_{\log^*}$.
  \item $\mathrm{PA} \vdash \mathrm{AR}_{f_\alpha} \Leftrightarrow \alpha < \varepsilon_0$.
  \item $\mathrm{I}\Sigma_{d+1} \vdash \mathrm{AR}^{d+1}_{f^{d+1}_\alpha} \Leftrightarrow \alpha < \omega_{d+2}$.
  \end{enumerate}
Where $f^{d+1}_\alpha (i) = \sqrt[H^{-1}_\alpha(i)]{\log^d (i)}$ and $f_\alpha (i)=\log^{H^{-1}_\alpha (i)}(i)$. 
\end{theorem} 

\begin{theorem}\label{thm:transition:PH}
  ~{}
  \begin{enumerate}
  \item $\mathrm{I}\Sigma_{d+1} \nvdash \mathrm{PH}^{d+2}_{\sqrt[c]{\log^d}}$ for every $c$. 
  \item $\mathrm{PA} \nvdash \mathrm{PH}_{\log^n}$ for every $n$.
  \item  $\mathrm{I}\Sigma_{1} \vdash \mathrm{PH}^{d+2}_{\log^{d+2}}$. 
  \item $\mathrm{I}\Sigma_1 \vdash \mathrm{PH}_{\log^*}$.
  \item $\mathrm{PA} \vdash \mathrm{PH}_{f_\alpha} \Leftrightarrow \alpha < \varepsilon_0$.
  \item $\mathrm{I}\Sigma_{d+1} \vdash \mathrm{PH}^{d+2}_{f^{d+1}_\alpha} \Leftrightarrow \alpha < \omega_{d+2}$.
  \end{enumerate}
Where $f^{d+1}_\alpha (i) = \frac{\log^{d+1} (i)}{H^{-1}_\alpha(i)}$ and $f_\alpha (i)=\log^{H^{-1}_\alpha (i)}(i)$. 
\end{theorem} 

\begin{theorem}\label{thm:transition:KM}
  ~{}
  \begin{enumerate}
  \item $\mathrm{I}\Sigma_{d+1} \nvdash \mathrm{KM}^{d+2}_{\sqrt[c]{\log^d}}$ for every $c$. 
  \item $\mathrm{PA} \nvdash \mathrm{KM}_{\log^n}$ for every $n$.
  \item  $\mathrm{I}\Sigma_{1} \vdash \mathrm{KM}^{d+2}_{\log^{d+1}}$. 
  \item $\mathrm{I}\Sigma_1 \vdash \mathrm{KM}_{\log^*}$.
  \item $\mathrm{PA} \vdash \mathrm{KM}_{f_\alpha} \Leftrightarrow \alpha < \varepsilon_0$.
  \item $\mathrm{I}\Sigma_{d+1} \vdash \mathrm{KM}^{d+2}_{f^{d}_\alpha} \Leftrightarrow \alpha < \omega_{d+2}$.
  \end{enumerate}
Where $f^{d}_\alpha (i) = \sqrt[H^{-1}_\alpha(i)]{\log^d (i)}$ and $f_\alpha (i)=\log^{H^{-1}_\alpha (i)}(i)$. 
\end{theorem} 
The first two items and the unprovability parts of the last two items of these theorems will be treated in Section~\ref{section:lowerbounds}. The first item of Theorem~\ref{thm:transition:AR},~\ref{thm:transition:PH} or~\ref{thm:transition:KM}  is Theorem~\ref{thm:independence:AR}, ~\ref{thm:independence:PH},~\ref{thm:independence:KM} combined with Theorem~\ref{thm:lowerbound:parAR},~\ref{thm:lowerbound:parPH},~\ref{thm:lowerbound:parKM} respectively. The unprovability parts of items (5) and (6) are shown by combining these with Lemmas~\ref{lemma:prooftheoreticsharpening} and~\ref{lemma:modeltheoreticsharpening}. 

Items (3) and (4) and the provability parts of the last two items are shown in Section~\ref{section:upperbounds}. These are direct consequences of Lemmas~\ref{lemma:upperbounds} and~\ref{lemma:upperboundssharpening} and upper bound estimates from the literature. 

Theorem~\ref{thm:transition:PH} can already be found in \cite{weiermann2004}, Theorem~\ref{thm:transition:KM} can be found in \cite{lee} and \cite{CLW}.
\section{Lower bounds}~\label{section:lowerbounds}
In the following three subsections we show items (1) and (2) of Theorems~\ref{thm:transition:AR},~\ref{thm:transition:PH} and ~\ref{thm:transition:KM}. The underlying idea of the three proofs is to show for the appropriate parameter $f$ that  $\varphi_f \rightarrow \varphi_{\mathrm{id}}$ by compressing the colourings $C$ for $\varphi_{\mathrm{id}}$ using $f$ to obtain a colouring $D_1$. This causes the problem that if one obtains for such colourings an adjacent/homogeneous/$\min$-homogeneous set (by $\varphi_f$) this set needs not satisfy adja-cency/homogeniety/$min$-homogeneiety (to demonstrate $\varphi_{\mathrm{id}}$) because $f(x)=f(y)$ may be satisfied for some $x<y$. We will solve this by combining  $D$ with two colourings $D_2$ and $D_3$. 

The colouring, $D_2$ will have the property that for adjacent/homogeneous/$\min$-homoge-neous set $H$ either $f(x_d)=f(x_{d+1})$ or $f(x_1) < \dots < f(x_{d+1})$ for $x_1 < \dots < x_{d+1}$ in $H$. 

The other colouring, $D_3$,will ensure that in the case $f(x_d)=f(x_{d+1})$ the set $H$ cannot be  adjacent/homogeneous/$\min$-homogeneous. To obtain the appropriate colouring $D_3$ we use lower bound estimates for $\varphi_k$ with constant function $k$. 

\subsection{Adjacent Ramsey}
For determining the transitions, estimates on $i \mapsto \mathrm{AR}_i^d(r)$ play a central part. A variant of these functions has been examined extensively in \cite{friedman2010}. We use the following result:
\begin{lemma} \label{lemma:estimate:ARconstantneq}
For every $d,i,c$ there exists a colouring 
\[
C \colon \{ 0, \dots , \mathrm{2}_d (i^c) \}^{d+1} \rightarrow \{0, \dots , i\}^{32 \cdot d+ c}
\]
such that $C(x) \neq C(y)$ for all $x \neq y$.
\end{lemma}
We modify this colouring slightly:
\begin{lemma} \label{lemma:estimate:ARc}
For every $d,c,i$ there exists a colouring 
\[
C_{d,c,i} \colon \{ 0, \dots , \mathrm{2}_d (i^c) \}^{d+1} \rightarrow \{0, \dots , i\}^{64 \cdot d+ 2\cdot c}
\]
such that $C(x_1, \dots , x_d) \not \leq C(x_2, \dots , x_{d+1})$ for all $x_1 < \dots < x_{d+1} \leq \mathrm{2}_d (i^c)$.
\end{lemma}
\emph{Proof:} Take $C'$ from Lemma~\ref{lemma:estimate:ARconstantneq} and define:
\[
C(x) = (C'(x), i-C'(x)).
\]
\begin{flushright} $\Box$ \end{flushright}
With these estimates we can prove parts (1) and (2) of Theorem~\ref{thm:transition:AR}. 
\begin{theorem}\label{thm:lowerbound:parAR}
There exists a  primitive recursive function $h$ such that:
\[
\mathrm{AR}_{\sqrt[c]{\log^d}}^{d+1}(h(d,c,r)) \geq \mathrm{AR}_{\mathrm{id}}^{d+1} (r).
\]
\end{theorem}
\emph{Proof:} We claim that the inequality holds for $h(d,c,r)=r+65\cdot (d+1) +2\cdot c+3$. Given $\mathrm{id}$-limited colouring $C \colon \{0, \dots , R\}^{d+1} \rightarrow \mathbb{N}^r$. Take $f(x)=\sqrt[c]{\log^d}$ and $(w(i))_j=1$ if $i = j$, zero otherwise.   

Define the following colourings:

$D_1 (x) = C(f(x_1) , \dots , f(x_{d+1})).$ 

$D_2 (x)=w(i)$,  where $i$ is the largest such that $f(x_{i-1})=f(x_i)$ if such $i$ exists,  one otherwise.

$D_3 (x)= C_{d+1, c+1, f(\max x)} (x)$, where $C_{d,c,i}$ are taken from Lemma~\ref{lemma:estimate:ARc}.

Combine these colourings into a single $f$-limited colouring
\[
D\colon \{ 0, \dots , R \}^{d+1} \rightarrow \mathbb{N}^{r+d+1+64\cdot (d+1) +2(c+1)}.
\] 
Suppose that for $x_1 < \dots < x_{d+2}$ we have $D(x_1, \dots , x_{d+1}) \leq D(x_2, \dots , x_{d+2})$. We have the following cases:
  \begin{enumerate}
  \item If $f(x_1) < \dots < f(x_{d+2})$ then, by definition of $D_1$, we have 
  \[
  C(f(x_1), \dots , f(x_{d+1})) \leq C(f(x_2), \dots , f(x_{d+2})).
  \]
  \item If there exists $i$ such that $f(x_{i-1}) = f(x_i) < \dots < f(x_{d+2})$ and $i=d+1$ then $D(x_1 , \dots , x_{d+1}) \not \leq D(x_2, \dots , x_{d+2})$ is inherited from $D_3$.
  \item If there exists $i$ such that $f(x_{i-1}) = f(x_i) < \dots < f(x_{d+2})$ and $i<d+1$ then $D(x_1 , \dots , x_{d+1}) \not \leq D(x_2, \dots , x_{d+2})$ is inherited from $D_2$.
  \end{enumerate}

\begin{flushright} $\Box$ \end{flushright}

\subsection{Paris--Harrington}
We will use lower bounds from Ramsey theory from  \cite{grarospe}, which are attributed to Erd\H{o}s and Hajnal:
\begin{lemma}\label{lemma:errad}
For every $d \geq 2$ there exists constant $a_d$ such that 
\[
\mathrm{PH}^{d}_{i \cdot a_d}(0,r) >\mathrm{2}_{d-2}(r^{i-2})
\]
for all $r \geq 4$ and $i \geq 3$. 
\end{lemma}
With these estimates we can prove parts (1) and (2) of Theorem~\ref{thm:transition:PH}. 
\begin{theorem}\label{thm:lowerbound:parPH}
There exist primitive recursive functions $h_1$ and $h_2$ such that:
\[
\mathrm{PH}^{d+1}_{\frac{\log^d}{c}} (h_1(c,d,m) , h_2(c,d,m,r))\geq \mathrm{PH}^{d+1}_{\mathrm{id}}(m,r)
\]
for every $c>0$ and $m$ sufficiently large.
\end{theorem}
\emph{Proof:}
We claim this is the case for: $h_1(c,d,m)= 2_d(a_d \cdot c \cdot m) $, $h_2(c,d,m,r)= r \cdot (d+2)^2 \cdot  2^{(c+2) \cdot a_d}$.

Given $C\colon [m,R]^{d+1} \rightarrow r$. Take $f(i)= \frac{\log^d(i)}{a_d \cdot c}$ and colourings , 
\[
D_{i\cdot a_d}\colon [\mathrm{2}_{d-1}(2^{(c+1) \cdot a_d \cdot i})]^{d+1} \rightarrow 2^{(c+2) \cdot a_d},
\]
where $D_{i\cdot a_d}$ is obtained from Lemma~\ref{lemma:errad}. Define
\[
D\colon  [2_d(a_d \cdot c \cdot m),R]^{d+1} \rightarrow r \times  2^{(c+2) \cdot a_d}\times (d+2)^2
\]
as follows: 

If $f(x_1) < \dots < f(x_{d+1})$ then:

$D(x)=(C( f(x_1) , \dots , f(x_{d+1})), 0, 0 , d+1)$.

If $1< i < d+1$ is the biggest $i$ such that $f(x_1) = \dots = f(x_i)$  and  $1 \leq  j < d+1$ is the biggest $j$ such that $f(x_1) < \dots < f(x_j)$, then:

$D(x)=(0,0, i, j)$.

Otherwise:

$D(x)=(0, D_{f(x_1) \cdot a_d}(x), d+1, 0)$.

Suppose that $H$ is homogeneous for $D$ and of size greater than $d+2$. In this case the last two coordinates have value $0$ or $d+1$. If not then there exist $x_1 < \dots < x_{d+2}$ with $i=(D(x_2, \dots x_{d+2}))_3=(D(x_1, \dots , x_{d+2}))_3+1=i+1$, contradiction (same argument for 4th coordinate). If one of those two is $0$ then the other must be $d+1$, so either $f(x_1)= \dots =f(x_{d+1})$ for all $x_1 < \dots < x_{d+1}$ in $H$  or $f(x_1)< \dots <f(x_{d+1})$ for all $x_1 < \dots < x_{d+1}$ in $H$.  

By definition of $D$ this implies that $H$ is homogeneous for $(D)_1$ or $D_{f(\min H)\cdot a_d}$. In the latter case $H$ has size strictly less than $f(\min H)$. Hence if $H$ has size larger than $f(\min H)$ then $H'=\{ f(h) : h \in H\} $ has size larger than $\min H'$ and is homogeneous for $C$.
\begin{flushright} $\Box$ \end{flushright}

\subsection{Kanamori--McAloon}
We have the following estimates from \cite{CLW}:
\begin{lemma} \label{lemma:estimate:KMc}
For every $d\geq 2$ there exists constant $a_d$ such that 
\[
\mathrm{KM}^d_{i \cdot a_d \cdot m }(0, a_d \cdot (m+1)) > \mathrm{2}_{d-2} (i^m).
\] 
\end{lemma} 
This implies that for $i> (a_d \cdot m)^m$ and $m>c+2$:
\[
\mathrm{KM}^d_i (0, a_d \cdot (m+1)) > \mathrm{2}_{d-2} ((i+1)^{c+1}).
\]
We denote bad colourings that show this with $D_i$. With these estimates we can prove parts (1) and (2) of Theorem~\ref{thm:transition:KM}. 

\begin{theorem}\label{thm:lowerbound:parKM}
There exist primitive recursive $h_1,h_2$ such that for $d \geq 2$:
\[
\mathrm{KM}^{d}_{\sqrt[c]{\log^{d-2}}}( h_1(d,m,c) , h_2(d,m,c) ) \geq \mathrm{KM}^{d}_{\mathrm{id}}(0,m).
\]
\end{theorem}
\emph{Proof:} We claim this inequality  holds for $h_1=\mathrm{2}_{d-2}((a_d \cdot m)^{cm})$, $h_2=a_d \cdot (m+1)+1$ and $m>d+c+3$.

Given a colouring $D\colon [R]^d \rightarrow \mathbb{N}$ for the identity function we create an intermediate colouring 
\[
\tilde{C}\colon  [R]^k \rightarrow \mathbb{N} \times \mathbb{N} \times (d+2) \times (d+2).
\]
Roughly speaking $\tilde{C}_1$ will be $D( f(x_1) , \dots , f(x_d))$,  $\tilde{C}_2$ is $D_{f(x_1)}$  and $\tilde{C}_{3,4}$ will ensure that for min-homogeneous sets either $f(x_1)= \dots =f(x_d)$ or $f(x_1)< \dots <f(x_d)$ in the manner similar to what we have seen for adjacent Ramsey and Paris--Harrington. We define $f$-regressive $C$ to be one of the first two coordinates or zero, where the choice is dependent on and coded by the value of the last coordinate. We emphasise again that the lower bound estimates for  $\mathrm{KM}^d_i$ directly influence the functions $f$ for which this construction is useful. 

We take $f=\sqrt[c+1]{\log^{d-2}}$ and:
\[
\tilde{C}(x)=(D( f(x_1) , \dots , f(x_d)), D_{f(x_1)}(x), i,j ), 
\]
where $i$ is the biggest such that $f(x_1)= \dots = f(x_i)$ and $j$ is the biggest such that $f(x_1)< \dots < f(x_j)$ ($j=1$ otherwise). Note that $\tilde{C}_1$ is not everywhere-defined, take  it to be $0$ if it is undefined (same for $\tilde{C}_2$). 

If $H$ of size at least $d+2$ is $\min$-homogeneous for $\tilde{C}_3$ then the values of this coordinate is $1$ or $d+1$. Suppose not,  let $x_1 < \dots <x_{d+1}$ be the first $d+1$ elements of $H$, then:
\[
i=\tilde{C}_3(x_1, x_3, \dots x_{d+1})=\tilde{C}_3(x_1, x_2, \dots x_d)+1=i+1,
\]
contradiction. 

If $H$ is $\min$-homogeneous for $\tilde{C}_4$ then it must by similar argument have values $1$, $2$ or $d+1$.  Let $x_1 < \dots < x_{d+1}$ be the first $d+1$ elements of $H$, and suppose $\tilde{C}_4(x_1, \dots , x_d)=2$, then $f(x_2)=f(x_3)$, hence  $\tilde{C}_4(x_2, x_3, \dots , x_d)=1$, in other words in this case $\tilde{C}_4$ has value $1$ on $H'=H-\min H$. 

Hence either $f(x) < f(y)$ for all $x<y \in H'$ or $f(x) = f(y)$ for all $x<y \in H'$. So $H'$ is $\min$-homogeneous for $D$ in the first case,  or $\min$-homogeneous for $D_{\min H'}$ in the latter case. 

Encode the last two coordinates into single colouring $E\colon [R]^d \rightarrow (d+1)^2$ such that the first of those two cases is encoded in value $0$, the latter in $1$. We take:
\[
C(x)=\left\{
\begin{array}{ll}
(d+1)^2 + 2 \cdot \tilde{C}_1(x)+1 & \textrm{if $E(x)=0$} \\
(d+1)^2 + 2 \cdot \tilde{C}_2(x)+2 & \textrm{if $E(x)=1$} \\
E(x)						  & \textrm{otherwise}.
\end{array}
\right.
\]
Suppose $H$ of size greater than $d+2$ is $\min$-homogeneous for $C$,  it must have value greater than $(d+1)^2+1$. Hence $H'=H-\min H$ is $\min$-homogeneous for either $D$ or $D_{f(\min H')}$. In the latter case it has size strictly less than $a_d\cdot(m+1)$.  Hence if we have a $\mathrm{min}$-homogeneous set for $C$ of size $a_d\cdot(m+1)+1$ we obtain a $\mathrm{min}$-homogeneous set for $D$ of size $m$.

This colouring is $\sqrt[c]{\log^{d-2}}$-regressive because 
\[
(d+1)^2+2+2 \cdot \sqrt[c+1]{\log^{d-2}(x_1)}< \sqrt[c]{\log^{d-2}(x_1)},
\]
is ensured by limiting the domain of $C$ to numbers larger than  $\mathrm{2}_{d-2}((a_d \cdot m)^{cm})$. 
\begin{flushright} $\Box$ \end{flushright}

\subsection{Sharpening}
In this subsection we prove the unprovability parts of (5) and (6) of Theorems~\ref{thm:transition:AR},~\ref{thm:transition:PH},~\ref{thm:transition:KM}. For applying the sharpening lemmas it is of use to note that if we combine the previous section with the results from \cite{friedmanpelupessy} and \cite{kanamorimcaloon} we have:
\begin{theorem}
Fix $d$. There exist primitive recursive $h_1, h_2, h_3, h_4,h_5$ such that
  \begin{enumerate}
  \item $M_{1,l_c^{-1}} (n,c,x) \geq H_{\omega_{d+1}(n)} (x)$ for 

  $M_{1,f} (n,c,x)= \mathrm{AR}^{d+1}_f(h_1(n,c,x))$ and $l_c(i)=\mathrm{2}_d (i^c)$.
  \item $M_{2,l_c^{-1}} (n,c,x) \geq H_{\omega_{d+1}(n)} (x)$ for 

  $M_{2,f} (n,c,x)= \mathrm{PH}^{d+2}_f(h_2(d,n,c,x), h_3(d, n,c,x))$ and $l_c(i)=\mathrm{2}_{d+1} (i\cdot (c))$.
  \item In every nonstandard model $N$ of $\mathrm{I}\Sigma_1$ with $n \in N$ and nonstandard $c,x \in N$ there exists a model of $\mathrm{I}\Sigma_{d}$ below $M_{3, l^{-1}_c}(n,c,x)$ for 
  
  $M_{3, f}(n,c,x)= \mathrm{KM}^{d+2}_f(h_4(d,n,c,x), h_5 (d,n,c,x))$ and $l_c(i)=\mathrm{2}_d (i^c)$.
  \end{enumerate}
\end{theorem}

\begin{lemma}[Proof theoretic lower bounds sharpening] \label{lemma:prooftheoreticsharpening}
Suppose $T$ is a theory which includes $\mathrm{I}\Sigma_1$ and we have the following:
  \begin{enumerate}
  \item $(i,c) \rightarrow l_c(i)$ is nondecreasing and provably total in $T$.
  \item $f(i) \leq g(i)$ for all $i \leq M_g(n,c,x)$ implies $M_f(n,c,x) \leq M_g(n,c,x)$.
  \item Every provably total function of $T$ can be eventually bounded by $H_n$ for some $n$ and $H(i)=H_i(i)$.
  \item $M_{l_c^{-1}} (n,c,i) \geq H_n (i)$ for every $n$.
  \end{enumerate}
then:
\[
T \nvdash \forall n,c,x \exists y M_h (n,c,x)=y,
\]
where $h(i)=l_{H^{-1}(i)}^{-1} (i)$.
\end{lemma}
\emph{Proof:} We show:
\[
M=M_h(x,x,x) \geq H(x).
\]
Suppose, for a contradiction, that:
\[
M < H(x),
\]
then $H^{-1} (i) \leq x$ for all $i \leq M$, hence $h(i) \geq l_x^{-1} (i)$ for all $i \leq M$. Therefore:
\begin{eqnarray*}
M & \geq &M_{l^{-1}_x} (x,x,x) \\
     & \geq & H_x(x) = H(x), 
\end{eqnarray*}
which contradicts our assumption.
\begin{flushright} $\Box$ \end{flushright}
\begin{corollary}~{}
  \begin{enumerate}
  \item $\mathrm{PA} \nvdash \mathrm{AR}_f$, where $f(i)=\log^{H_{\varepsilon_0}^{-1}(i)}(i)$.
  \item $\mathrm{I}\Sigma_{d+1} \nvdash \mathrm{AR}_f^{d+1}$ where $f(i)=\sqrt[H^{-1}_{\omega_{d+1}}(i)]{\log^d(i)}$.
  \item $\mathrm{PA} \nvdash \mathrm{PH}_f$, where $f(i)=\log^{H_{\varepsilon_0}^{-1}(i)}(i)$.
   \item $\mathrm{I}\Sigma_{d+1} \nvdash \mathrm{PH}_f^{d+2}$ where $f(i)=\frac{\log^{d+1}(i)}{H^{-1}_{\omega_{d+2}}(i)}$.
  \end{enumerate}
\end{corollary}

\begin{lemma}[Model theoretic lower bounds sharpening] \label{lemma:modeltheoreticsharpening}
Suppose $T$ is a theory which includes $\mathrm{I}\Sigma_1$ and we have the following:
  \begin{enumerate}
  \item $(i,c) \rightarrow l_c(i)$ is nondecreasing and provably total in $T$.
  \item $f(i) \leq g(i)$ for all $i \leq M_g(n,c,x)$ implies $M_f(n,c,x) \leq M_g(n,c,x)$.
  \item $H$ eventually dominates every provably total function of $T$.
  \item In every nonstandard model $N$ of $\mathrm{I}\Sigma_1$ and for every $c \in N$ and nonstandard $n,x, M_{l_c^{-1}}(n,c,x) \in N$ there exists an initial segment $I<M_{l_c^{-1}}(n,c,x)$ which models $T$.
  \end{enumerate}
then:
\[
T \nvdash \forall n,c,x \exists y M_h (n,c,x)=y,
\]
where $h(i)=l_{H^{-1}(i)}^{-1} (i)$.
\end{lemma}
\emph{Proof:} Fix nonstandard model $N$. If $M_h(x,x,x) \geq H(x)$ for infinitely many standard $x$ we are finished, so suppose that for all but finitely many standard $x$ we have:
\[
M_h(x,x,x) < H(x).
\]
For these $x$ we know that $h(i) \geq l_x^{-1}(i)$ for all $i \leq M_h(x,x,x)$, by overflow there exists nonstandard $x$ with these properties, so there exists a nonstandard instance of $M_{l_x^{-1}} (x,x,x)$.  Hence there exists initial segment which models $T$.
\begin{flushright} $\Box$ \end{flushright}
\begin{corollary}~{}
  \begin{enumerate}
  \item $\mathrm{PA} \nvdash \mathrm{KM}_f$, where $f(i)=\log^{H_{\varepsilon_0}^{-1}(i)}(i)$.
  \item $\mathrm{I}\Sigma_{d+1} \nvdash \mathrm{KM}_f^{d+2}$ where $f(i)=\sqrt[H^{-1}_{\omega_{d+2}}(i)]{\log^d(i)}$.
  \end{enumerate}
\end{corollary}

\section{Upper bounds}~\label{section:upperbounds}
In this section we show items (3) and (4) and the provability parts of items (5) and (6) of Theorems~\ref{thm:transition:AR},~\ref{thm:transition:PH} and ~\ref{thm:transition:KM}.
\begin{lemma}[Upper bounds lemma] \label{lemma:upperbounds}
Suppose $T$ is a theory that contains $\mathrm{I}\Sigma_1$, $M_f\colon\mathbb{N}^2 \rightarrow \mathbb{N}$ is a computable function for all computable $f$ and $M_f(d,x) \leq M_g(d,x)$ whenever $f(i)\leq g(i)$ for all $i \leq M_g(d,x)$. Additionally, suppose that there exist increasing, provably total, functions $u$,$h$ such that for every $d,n$ and $k \geq h(d,n)$ we have:
\[
M_k(d,n) \leq u(k),
\] 
then: 
\begin{eqnarray*}
T & \vdash 		&\forall d,x \exists !y M_{u^{-1}} (d,x)=y.
\end{eqnarray*}  
\end{lemma}
\emph{Proof:} If $i \leq u(h(d,x))$ then $u^{-1}(i) \leq h(d,x)$. Hence:
\[
M_{u^{-1}} (d,x) \leq M_{h(d,x)} (d,x) \leq u(h(d,x)).
\]
\begin{flushright} $\Box$ \end{flushright}
\begin{corollary}
If $\varphi$ is one of $\mathrm{AR}$, $\mathrm{PH}$, $\mathrm{KM}$ then:
\[
\mathrm{I}\Sigma_1 \vdash \varphi_{\log^*}
\]
\end{corollary}
\emph{Proof:} Ramsey numbers are bound by the tower function by the Erd\H{o}s-Rado bounds from \cite{erdosrado1951} .
\begin{flushright} $\Box$ \end{flushright}
\begin{lemma}[Upper bounds sharpening lemma] \label{lemma:upperboundssharpening}
Let $T,M$ be as in the upper bounds lemma. If $(c,i) \mapsto l_c(i)$ is an increasing provably total function such that there exist provably total $g_1,g_2$ with $g_1(d) \leq g_2(d,x)$ for all $x$ and $M_k(d,x)\leq l_{g_1(d)}(k)$ whenever $k\geq g_2(d,x)$, then:
\begin{eqnarray*}
T & \vdash 		&\forall d,x \exists !y M_{f} (d,x)=y,
\end{eqnarray*}  
where $f(i)=l_{B^{-1}(i)}^{-1}(i)$ and $B$ is an arbitrary unbounded, increasing and provably total function.
\end{lemma}
\emph{Proof:} If $i \leq l_{g_1(d)} (B(g_2(d,x)))$ then:
\[
f(i) \leq l^{-1}_{g_2(d,x)}(l_{g_1(d)} (B(g_2(d,x)))) \leq l^{-1}_{g_2(d,x)}(l_{g_2(d,x)}(B(g_2(d,x)))).
\]
Therefore:
\[
M_f(d,x) \leq M_{B(g_2(d,x))} (d,x) \leq  l_{g_1(d)} (B(g_2(d,x))).
\]
\begin{flushright} $\Box$ \end{flushright}

\begin{corollary}
$\mathrm{I}\Sigma_{d+1} \vdash \mathrm{AR}^{d+1}_{f_\alpha}$ whenever $f_{\alpha} (i)=\sqrt[H^{-1}_{\alpha}(i)]{\log^d(i)}$ and $\alpha < \omega_{d+2}$.
\end{corollary}
\emph{Proof:} Examine $\mathrm{AR}^d_f$ with fixed $d$ and its associated function $\mathrm{AR}^d_f(r)$. The $r$ will have the role of $d$ when applying the upper bounds sharpening lemma. By the Erd\H{o}s-Rado bounds on Ramsey numbers from \cite{erdosrado1951}:
\[
\mathrm{AR}_k^{d+1} (r) \leq \mathrm{2}_d (k^{(r+1)}).
\]
Hence, by sharpening, $\mathrm{I}\Sigma_{d+1} \vdash \mathrm{AR}^{d+1}_{f_\alpha}$ whenever $f_{\alpha} (i)=\sqrt[H^{-1}_{\alpha}(i)]{\log^d(i)}$ and $\alpha < \omega_{d+2}$.
\begin{flushright} $\Box$ \end{flushright}
\begin{corollary}
$\mathrm{I}\Sigma_{d} \vdash \mathrm{PH}^{d+1}_{f_\alpha}$ whenever $f_{\alpha} (i)=\frac{\log^d(i)}{H^{-1}_{\alpha}(i)}$ and $\alpha < \omega_{d+1}$.
\end{corollary}
\emph{Proof:}
Examine $\mathrm{PH}^d_f$ with fixed $d$ and its associated function $\mathrm{PH}^d_f(m,r)$. The $r$ will have the role of $d$ when applying the upper bounds sharpening lemma. By the Erd\H{o}s-Rado bounds on Ramsey numbers from \cite{erdosrado1951} if $k \geq r+m$ then:
 \[
\mathrm{PH}^d_k(m,r) \leq \mathrm{2}_{d-1} (r^{d^2} \cdot k)+m \leq \mathrm{2}_{d-1} ((r^{d^2} +1)\cdot k)=l_{r^{d^2}+1}(k).
\]
Hence, by sharpening, $\mathrm{I}\Sigma_{d} \vdash \mathrm{PH}^{d+1}_{f_\alpha}$ whenever $f_{\alpha} (i)=\frac{\log^d(i)}{H^{-1}_{\alpha}(i)}$ and $\alpha < \omega_{d+1}$.
\begin{flushright} $\Box$ \end{flushright}

\begin{corollary}
$\mathrm{I}\Sigma_{d+1} \vdash \mathrm{KM}^{d+2}_{f_\alpha}$ whenever $f_{\alpha} (i)=\sqrt[H^{-1}_{\alpha}(i)]{\log^d(i)}$ and $\alpha < \omega_{d+2}$.
\end{corollary}
\emph{Proof:}
Examine $\mathrm{KM}^d_f$ with fixed $d$ and its associated function $\mathrm{KM}^d_f(a,m)$.  The $m$ will have the role of $d$ when applying the upper bounds sharpening lemma.  We use bounds from Corollary 4.2.3 in \cite{lee}:
\[
\mathrm{KM}^d_k(a,m) \leq \mathrm{2}_{d-2} (k^{d^2\cdot m})+a \leq  \mathrm{2}_{d-2} (k^{d^2 \cdot m+2})=l_{d^2 \cdot n+2} (k),
\]
where the second inequality is true for $k \geq (a+d^2 \cdot m +2)$. Hence, by sharpening, $\mathrm{I}\Sigma_{d+1} \vdash \mathrm{KM}^{d+2}_{f_\alpha}$ whenever $f_{\alpha} (i)=\sqrt[H^{-1}_{\alpha}(i)]{\log^d(i)}$ and $\alpha < \omega_{d+2}$.
\begin{flushright} $\Box$ \end{flushright}
\begin{corollary}
Let $\varphi$ be one of $\mathrm{AR}$, $\mathrm{PH}$, $\mathrm{KM}$ and $f_\alpha(i)= \log^{H^{-1}_\alpha(i)}(i)$. We have:
\[
\mathrm{PA} \vdash \varphi_{f_\alpha},
\]
whenever $\alpha < \varepsilon_0$.
\end{corollary}
\section{Some observations on transitions} \label{section:observation}
In the phase transitions which have been examined so far the same heuristics are used to determine the threshold functions: as soon as the upper bound lemmas cannot be applied because $l$ is a lower bound the resulting theorem is not provable for $l^{-1}$. We conjecture that phase transitions in unprovability always have this shape:
\begin{conjecture}[Lower bounds]
Suppose $T$ is a theory which contains $\mathrm{I}\Sigma_1$, $l$ is nondecreasing and $M_f$ is a computable function for every computable $f$ with the following properties:
  \begin{enumerate}
  \item $T \nvdash \forall x \exists y M_{\mathrm{id}}(x)=y$,
  \item $f(i) \leq g(i)$ for all $i \leq M_g(x)$ implies $M_f(x) \leq M_g(x)$,
  \item There exists $x$ such that $k \mapsto l(k)$ is eventually bounded by $k \mapsto M_k (x)$, 
  \end{enumerate}
then:
\[
T \nvdash \forall x \exists y M_{l^{-1}}(x)=y.
\]
\end{conjecture}
For the sharpening of the transition results the two lower bounds sharpening lemmas suffice. These lemmas are dependant on the method of proving independence of $\varphi_{\mathrm{id}}$. We conjecture that it is possible to do this independent of this method: 
\begin{conjecture}[Lower bounds sharpening]
Suppose $T$ is a theory which contains $\mathrm{I}\Sigma_1$, $(c,i) \mapsto l_c(i)$ is nondecreasing and $M_f$ is a computable function for every computable $f$ with the following properties:
  \begin{enumerate}
  \item $T \nvdash \forall x \exists y M_{l^{-1}_c}(x)=y$ for every $c$,
  \item $f(i) \leq g(i)$ for all $i \leq M_g(x)$ implies $M_f(x) \leq M_g(x)$,
  \item $H$ eventually dominates every provably total function of $T$.
  \end{enumerate}
then:
\[
T \nvdash \forall x \exists y M_h (x)=y,
\]
where $h(i)=l_{H^{-1}(i)}^{-1} (i)$.\end{conjecture}

\end{document}